\documentclass[12pt,a4paper]{amsart}

\usepackage[utf8]{inputenc}
\usepackage{geometry}
\usepackage{amsmath}
\usepackage{amsfonts}
\usepackage{amssymb}
\usepackage{graphicx}
\usepackage{amsthm}
\usepackage{mathtools}
\usepackage{array}
\usepackage{enumerate}
\usepackage{float}
\usepackage{tikz-cd}
\usepackage[color = blue!20]{todonotes}
\usepackage{hyperref}
\usepackage[capitalise,noabbrev]{cleveref}
\usepackage{array,multirow}
\usepackage{pdflscape}

\numberwithin{equation}{section}
\theoremstyle{definition}
\newtheorem{theorem}{Theorem}[section]

\DeclareMathOperator{\SL}{\mathrm{SL}}

\newcommand{\Z}{\mathbb{Z}}

\newcommand{\Q}{\mathbb{Q}}

\newcommand{\Mod}[1]{\ (\mathrm{mod}\ #1)}
\newcommand{\demph}[1]{\textbf{#1}}

\newcommand{\abs}[1]{\lvert #1 \rvert}
\newcommand{\HH}{\mathcal{H}}

\renewcommand{\i}{\mathrm{i}}

\renewcommand\arraystretch{0.667}

\newcommand{\LMRefDB}[2]{\href{#1}{\texttt{#2}}}
\title{Expressions for weight 2 cusp forms in holomorphic eta quotients}

\begin{document}

\author[Chan]{Elisabeth (Yin Ting) Chan}
\address{Elisabeth (Yin Ting) Chan}
\email{\url{elisabeth15501@gmail.com}}

\author[Combes]{Lewis Combes}
\address{Lewis Combes, School of Mathematics and Statistics, University of Sheffield, UK, S3 7RH}
\email{\url{lmc577@proton.me}}

\begin{abstract} 
We attempt to compute expressions in terms of the Dedekind eta function for all weight 2 new cusp forms with level up to 100, using methods of Allen et. al. In cases where no expression exists, we raise the level instead of the weight, meaning our eta quotients are always holomorphic. Of the forms we examine, we find expressions for all but 4. We also present methods to find expressions with relatively few terms, and how these expressions can be used to demonstrate zeroes of modular forms.
\end{abstract}

\maketitle
\tableofcontents

\section{Introduction}
Let $\HH_2 = \left\{ x+y\i  \mid y > 0 \right\}$ be the complex upper half-plane. We will write $M_k(N)$, $S_k(N)$ for the spaces of modular and cuspidal forms of weight $k$ and level $\Gamma_0(N)$ on $\HH_2$. The \demph{Dedekind eta function} is defined for $z \in \HH_2$ as 
\begin{equation}\label{EqnEtaExp}
    \eta(z) = q^{1/24}\prod_{n=1}^{\infty} (1-q^n),
\end{equation}
where $q=e^{2\pi \i z}$. It is a modular form for the full modular group $\Gamma = \SL_2(\Z)$ of weight $1/2$ (see \cref{SectionZeroesOfMFs}). Many modular forms can be expressed in terms of the eta function, for example: 
\begin{equation*}
    \Delta(z) = \eta(z)^{24}, \hspace{1pc} f(z) = \eta(z)^2\eta(11z)^2, \hspace{1pc} g(z) = \eta(z)^8\eta(4z)^4 + 8\eta(4z)^{12}.
\end{equation*}
Here $\Delta$ is the discriminant cusp form of weight $12$ and level $1$; $f$ is the unique cusp form of weight $2$ and level $11$; $g$ is the unique newform of weight $6$ and level $8$. These forms have LMFDB labels \LMRefDB{https://www.lmfdb.org/ModularForm/GL2/Q/holomorphic/1/12/a/a/}{1.12.a.a}, \LMRefDB{https://www.lmfdb.org/ModularForm/GL2/Q/holomorphic/11/2/a/a/}{11.2.a.a} and \LMRefDB{https://www.lmfdb.org/ModularForm/GL2/Q/holomorphic/8/6/a/a/}{8.6.a.a} respectively (see \cite{lmfdb:cmf.label} for details on these labels). The first of these two expressions are well-known, and follow from the formulae in \cref{SectionImportantResults} and the fact that $S_{12}(1)$ and $S_{2}(11)$ are 1-dimensional. The third can be found in \cite{OsburnStraubZudilin}. 

An expression of the form $\prod_{d\mid N} \eta(dz)^{r_d}$ is called an \demph{eta quotient}. Eta quotients crop up often in the study of modular forms, partly due to their utility in specifying modular forms concisely. Instead of referring to ``the unique cuspidal newform of level 8 and weight 4'', one can simply write $\eta(2z)^4 \eta(4z)^4$ (a fact proved in \cite{MartinEtaQuots}). 

We are motivated by this aim of using eta quotients to find concise ways to write down modular forms. In particular, we aim to find expressions for weight 2 forms associated to elliptic curves (i.e. the forms that are both new at their level and cuspidal) in terms of the eta function. This particular class of forms is chosen for the relatively low dimension of their spaces and their ease of computation. 

Expressing \textit{all} such forms as eta quotients is immediately hopeless, since only finitely many of these forms have such expressions---a fact proved by Martin and Ono \cite{MartinOno}. In general, then, one can only hope to express a given form as a \textit{linear combination} of eta quotients, which we call an \demph{eta expression}. 

After this relaxation, many more forms can be expressed using the eta function. For example, Allen et. al. \cite{AllenEtAl2020} give an eta expression for the unique cusp form of weight 2 and level 35:
\begin{equation*}
     f_{35}(z) = \eta(z)^2\eta(35z)^2 + \eta(5z)^2\eta(7z)^2,
\end{equation*}
as well as a longer expression of the form \LMRefDB{https://www.lmfdb.org/ModularForm/GL2/Q/holomorphic/55/2/a/a/}{55.2.a.a} at level 55, and a method to generate more expressions. 

A modular form of a given weight $k$ and level $N$ need not have an expression in terms of eta quotients living in $M_k(N)$; in particular, if $N$ has few factors, the number of eta quotients in $M_k(N)$ may not even span the space. The expression of \cite{AllenEtAl2020} bypasses this problem by raising the weight, allowing spaces with more eta quotients to be considered. We briefly recap the idea here. 

If $f \in M_k(N)$ is not expressible as a linear combination of eta quotients in $M_k(N)$, we can pick some eta quotient $g \in M_l(N)$ and consider the product $fg$, which may be expressible as a linear combination $l = \sum a_i \eta_i$ of eta quotients in $M_{k+l}(N)$. Then $l/g = \sum a_i \eta_i/g$ is also a linear combination of eta quotients, representing $f$. However, each of the $\eta_i/g$ need not be holomorphic at all the cusps of $X_0(N)$, and so may not be in $M_{k}(N)$; indeed, if $f$ \textit{can't} be expressed as a linear combination of eta quotients at its own weight and level, then at least one of the quotients cannot be holomorphic at all the cusps. Such quotients are called \demph{weakly-holomorphic}. 

In this paper, we will take a complementary approach: if $f$ is not in the span of the eta quotients at its own level $N$, we instead treat it as an element of $M_2(Nd)$ for some $d>1$, chosen so that $f$ is in the span of the eta quotients of level $Nd$. In this way, all of our eta quotients remain holomorphic at the cusps, at the cost of raising the level. This seems to have a benefit over the weight-raising approach, as the expression we find for \LMRefDB{https://www.lmfdb.org/ModularForm/GL2/Q/holomorphic/55/2/a/a/}{55.2.a.a} is made of four eta quotients with small coefficients, at the cost of an extra factor of $2$ in the level; meanwhile, the expression found by raising the weight in \cite{AllenEtAl2020} is made of 39 eta quotients, with complicated fractions for coefficients. 

In this way, we are particularly interested in finding expressions that are \demph{human-readable} (to the extent that this is possible). We illustrate with an example. Let $f$ be the form with label \LMRefDB{https://www.lmfdb.org/ModularForm/GL2/Q/holomorphic/32/2/a/a/}{32.2.a.a}. The space $M_2(32)$ has dimension 8, and contains 131 eta quotients, meaning there is a great deal of linear dependence between the quotients. From \cite{MartinOno}, we know that 
\begin{equation*}
    f(z) = \eta(4z)^2\eta(8z)^2,
\end{equation*}
but it can also be computed, using the methods of \cref{SectionFindingExpressions} that 
\begin{equation*}
    f(z) = \frac{ \eta(z)^2 \eta(4z)^2 \eta(16z)^2 }{ \eta(2z) \eta(8z)^3 } + 2 \frac{\eta(4z)^4}{\eta(8z)^4}. 
\end{equation*}
Clearly the first expression is the more human-readable of the two. We formalise this notion via \textit{minimality} in \cref{SubsecMinimalExpressions}, which guides our search for eta expressions. Very often it is easy to raise the level of a form until it is a linear combination of some eta quotients, but there is no guarantee that this expression is minimal; generally, the first expression found is long and unreadable, and requires work to refine. We describe algorithms to do this in \cref{SubsecFindingMinimalExpressions}.

\subsection{Summary of computations}

The aim of the original project from which this paper followed was to find expressions for all the cusp forms associated to elliptic curves of conductor up to 100, as a modest test case for more general computations. There are 93 of these forms, some of which have expressions floating around in the literature. 

Collecting all that we can find, we have 19 eta expressions from various sources that are minimal in our sense\footnote{There are a few more that can be found in \cite{AyginThesis}, but we find these are not minimal.}. Of the remaining 74, we find expressions for 70, as well as bounds on the minimal factor needed to raise the level by to get an expression for the other 4. All the expressions that are minimal (as well as some non-minimal expressions of reasonable length) can be found in the tables of \cref{SectionTables}. Of the 70 new expressions we find, 58 are minimal. The code used to compute all expressions, eta quotients for all the levels considered here, and some larger expressions, can all be found at the GitHub repository
\begin{center}
    \url{https://github.com/lewismcombes/EtaExpressions}
\end{center}

\subsection{Structure of the paper}

In \cref{SectionImportantResults}, we collect some important results in the study of eta quotients, used to generate complete lists of all the eta quotients of a given weight and level. In \cref{SectionFindingExpressions}, we describe algorithms to find eta expressions, both generally and with the restriction of being minimal. In \cref{SectionZeroesOfMFs}, we use one of our eta expressions to prove the location of a zero of a modular form, in a fairly mechanical process that could likely be refined to provide many more zeroes from eta expressions. In \cref{SectionTables}, we list the eta expressions we have found.

\subsection{Acknowledgements} 
This paper grew out of the University of Sheffield's Undergraduate Research Internship project \textit{``Eta expressions associated to elliptic curves''}, completed during summer 2022 by the first author and supervised by the second. The UGRI scheme was funded by a grant from the Heilbronn Institute for Mathematical Research, and the second author was funded as a PhD student during the work by the EPSRC grant EP/R513313/1. The computations were performed in \textsf{Magma} V2.28-5, on the Sheffield School of Mathematics and Statistics server, Intel\textsuperscript{\textregistered} Xeon\textsuperscript{\textregistered} CPU E5-2683 v3 @ 2.00GHz, with 128GB of RAM, running Ubuntu 22.04.2 LTS.

\section{Important results}\label{SectionImportantResults}

We use the standard convention to write eta quotients compactly. Let $d_1, \ldots, d_l$ be the divisors of $N$, ordered by size. We write 
\begin{equation*}
    \eta_{N}[r_1,\ldots,r_l] = \prod_{i=1}^l \eta(d_i z)^{r_i}. 
\end{equation*}
Not every eta quotient of level $N$ is a modular form for $\Gamma_0(N)$. Those that are can be classified using the following results. First, Newmann \cite{Newman2} gives the following. 
\begin{theorem}[Newmann]
    Let $f(z) = \eta_N[r_1, \ldots, r_l]$ be an eta quotient, with 
    \begin{equation}\label{EqnDivisorSumCond1and2}
        \sum_{i=1}^{l} d_i r_i \equiv 0 \Mod{24}, \hspace{1pc}  \sum_{i=1}^{l} \frac{N}{d_i}r_i \equiv 0 \Mod{24},
    \end{equation}
    and
    \begin{equation}\label{EqnProdCond}
        \prod_{i=1}^{l} d_i^{r_i} \text{ is a square in } \Q.
    \end{equation}
    Then $f(z)$ transforms like a modular form of weight $k= \frac{1}{2} \sum  r_i $ for $\Gamma_0(N)$. 
\end{theorem}
Any eta quotient satisfying these three conditions can fail to be a modular form exactly when it is not ``holomorphic at the cusps'' of $\Gamma_0(N)$, i.e. if $f(z)$ is not bounded as $z$ approaches a $\Gamma_0(N)$-equivalence class of points in $\Q \cup \{ \infty \}$. This can be measured with the following result of Ligozat \cite{Ligozat}. 
\begin{theorem}[Ligozat]
    Let $f(z) = \eta_N[r_1, \ldots, r_l]$ be an eta quotient. The order of vanishing of $f(z)$ at the cusp $s = \frac{c}{d}$ is given by 
    \begin{equation*}
        v_{s}(f) = \frac{N}{24}\sum_{i=1}^l \frac{\gcd(d,d_i)^2r_i}{\gcd(d,N/d)\cdot d \cdot d_i}. 
    \end{equation*}
\end{theorem}
These results combine to give an effective process for listing all the eta quotients of a given level and weight. This has been implemented by Rouse and Webb \cite{RouseWebb}, in particular in their \textsf{Magma} script\footnote{Available at \url{https://users.wfu.edu/rouseja/eta/}} \texttt{etamake-lattice.txt}. It is this script we used to compute our lists of holomorphic eta quotients.

\section{Finding eta expressions}\label{SectionFindingExpressions}
Once a list of eta quotients has been found, it is relatively straightforward to express a chosen $f(z) \in M_k(N)$ as a linear combination. Using \cref{EqnEtaExp}, the $q$-expansion of each quotient can be obtained to arbitrary precision. These can then be identified as elements of $M_k(N)$ by ensuring the $q$-expansion is found up to the $B$\textsuperscript{th} term, where $B$ is the Sturm bound (a positive integer past which the $q$-expansions of different modular forms must differ). Once a representation of $f$ and each eta quotient is found in terms of some basis of $M_k(N)$, the rest is linear algebra. 

If one is only interested in the \textit{existence} of eta quotients, this method is sufficient, and the rest of this section can be skipped. However, this approach has the downside of often returning eta expressions that are the sum of many terms, which can be unwieldy to write down, and defeat the purpose of expressing modular forms succinctly. Therefore we are interested in finding expressions that are \demph{minimal}, in a sense we now describe. 

\subsection{Minimal expressions}\label{SubsecMinimalExpressions}
Suppose $f(z)$ has two expressions as a linear combination of eta quotients $g_i$, $h_j$: 
\begin{equation*}
    f(z) = \sum_{i=1}^n a_i g_i(z) = \sum_{j=1}^m b_j h_j(z). 
\end{equation*}
We say the first expression is smaller than the other if 
\begin{enumerate}
    \item $n < m$, or 
    \item $n = m$ and $\sum_{i=1}^n\mathcal{H}(a_i) < \sum_{j=1}^m \mathcal{H}(b_j)$,
\end{enumerate}
where $\mathcal{H}$ is the height function on rationals $x = \frac{a}{b}$ given by 
\begin{equation*}
    \mathcal{H}(x) = \log(\abs{ab}). 
\end{equation*}
We note, this is different from the usual notion of height given by $\max\{ \abs{a},\abs{b} \}$. We make this choice because we want to account for both the written length of the numerator and the denominator of the coefficients. 

A minimal eta expression is therefore one for which no smaller expression exists. It is important to note that the above notion of size only imposes a partial order on the set of eta expressions of $f$; minimal expressions need not be unique, but they must exist if $f$ has at least one eta expression. This does not necessarily mean that a minimal expression is overly ``nice'' to write down, for example the expression for \LMRefDB{https://www.lmfdb.org/ModularForm/GL2/Q/holomorphic/26/2/a/a/}{26.2.a.a} in \cref{SectionTables} is minimal, but is still not very pleasant to look at. 

This example raises one final important note on minimality: when $f$ is not expressible as an eta quotient at its own level $N$, we express it at a higher level $Nd$, choosing the minimal $d$ such that this can be done. It is conceivable that an expression for $f$ could be found at a higher level that uses fewer quotients than a minimal expression at a lower level; we will put this possibility to one side, as a possible consideration for future research. 

\subsection{Algorithms for minimal expressions}\label{SubsecFindingMinimalExpressions}
We use two main algorithms to bound the length of a minimal eta expression from above and below. In the best cases, these two bounds meet in the middle, and we can identify a minimal expression exactly. 

\subsubsection{Bounding below}
There is a simple algorithm for finding a minimal eta expression, although it is unsophisticated and extremely inefficient. Assuming one has found $f$ and every eta quotient as vectors in $M_k(N)$, one can check each subspace of $M_k(N)$ generated by 1 eta quotient to see if it contains $f$, then each space generated by 2 eta quotients, and so on. Carrying on in this way, the algorithm is guaranteed to find the minimal such $n$ such that $f$ is a linear combination of $n$ eta quotients; then all the expressions can be checked against each other to find one of minimal height. 

The obvious drawback of this algorithm is that it is exhaustive, and when there are many eta quotients in $M_k(N)$ it becomes infeasible in practice. If there are $m$ eta quotients in $M_k(N)$, there are $\binom{m}{n}$ subspaces generated by $n$ eta quotients, and even for modest values of $m$ and $n$, this number quickly gets out of hand. 

There are many levels in our tables where we \textit{can} guarantee the minimality of an expression using this algorithm; for those where $\binom{m}{n}$ becomes too large, we only guarantee that a minimal expression doesn't exist for $n$ in certain ranges. As a simple cutoff, we only check all possible combinations up to the point that $\binom{m}{n} \leq 10^8$. 

\subsubsection{Bounding above}
We use a random algorithm to bound the length of a minimal eta expression from above. The method is simple: we build a subspace of $M_k(N)$ one eta quotient at a time, chosen at random, checking each time whether it contains $f$. Once a space containing $f$ is found, we check it against the current record holder for minimality, discarding the record holder if a smaller expression is found. 

While fairly simple, this algorithm has proved very effective in practice, allowing expressions to be found in cases where the exhaustive method of bounding from below becomes ineffective. For example, the expressions for the forms of level 90 in \cref{SectionTables} were found using this method. While they cannot currently be guaranteed to be minimal, they are certainly close to being so.

\subsubsection{Inferred minimality}
In general, one needs to find every eta expression for $f$ of length $n$ to guarantee that a given expression of length $n$ is minimal, since the coefficients of a minimal expression must be checked against all others. However, in some special cases we can infer the minimality of an eta expression using its known bounds. If we have an eta expression of length $n$ whose coefficients are all $\pm 1$, and we have verified that no expression of length $n-1$ exists, then we know our expression must be minimal since no coefficients with smaller height can exist. There may be other expressions of length $n$ whose coefficients are \textit{also} $\pm 1$, but this doesn't matter to us, as long as we have found one. 

For example, at level $96$ there are 22,655 eta expressions, meaning one should naively have to check $2.6\times 10^{8}$ subspaces of dimension 2 to prove the expressions we give are minimal, which is outside of our computational cutoff. However, we know that no expression of length 1 can exist for either form (due to the completeness of the list of \cite{MartinOno}), and the coefficients in the expressions we find are all $\pm 1$, so we can certify the expression as minimal without checking all cases.

\subsection{A speculative heuristic}
Here we note a speculation used to find the expression we list for \LMRefDB{https://www.lmfdb.org/ModularForm/GL2/Q/holomorphic/72/2/a/a/}{72.2.a.a}. Looking through the tables, we find many levels where a minimal expression is given as the sum/difference of two eta quotients $e_1 = \eta_N[d_1, \ldots, d_l]$, $e_2 = \eta_N[\delta_1, \ldots, \delta_l]$, where the $d_i$ are a permutation of the $\delta_j$. We do not know of any reason this should happen as often as it does, and it is certainly not true that we \textit{always} find a minimal expression of this form.

Nevertheless, it turned out to work in this instance. By taking subsets of our list of eta expressions at weight $2$ and level $72$ whose entries are permutations of the same (multi-)set of numbers, we found the expression listed in a few minutes. The exhaustive algorithm on all 45,798 eta expressions would have required checking $\approx 10^9$ subspaces of dimension 2. 

We also note that some expressions are comprised of multiple pairs of permuted quotients, such as the expression we list for \LMRefDB{https://www.lmfdb.org/ModularForm/GL2/Q/holomorphic/39/2/a/a/}{39.2.a.a}. This is not always the case, but seems to happen more often than chance should allow for. In the case of this particular expression, the coefficients of the individual quotients in the ``permuted pairs'' are the same (up to sign). Again, this does not happen in general, but poses an intriguing question as to why this seems to happen more often than it should.

\section{Application: zeroes of modular forms}\label{SectionZeroesOfMFs}
Let $f$ be the cusp form with LMFDB label \LMRefDB{https://www.lmfdb.org/ModularForm/GL2/Q/holomorphic/42/2/a/a/}{42.2.a.a}. Per \cref{SectionTables}, we have 
\begin{equation*}
    f(z) = \eta_{42}[ -1, 2, 2, -1, -1, 2, 2, -1 ] - \eta_{42}[ 2, -1, -1, 2, 2, -1, -1, 2 ].
\end{equation*}
Using a root-finding algorithm\footnote{Such as, for example, \textsf{Mathematica}'s \texttt{FindRoot}} on this expression, it appears that $f(w)=0$, where
\begin{equation*}
    w = \frac{1}{2} + \frac{1}{2\sqrt{21}}\i. 
\end{equation*}
Manipulating the eta expression for $f$ above, we see that $f(w)=0$ if and only if 
\begin{equation}\label{Equation42EtaRelation}
    \eta(w)^3 \eta(6w)^3 \eta(7w)^3 \eta(42w)^3 = \eta(2w)^3\eta(3w)^3\eta(14w)^3\eta(21w)^3,
\end{equation}
since $\eta$ is non-zero on $\HH_2$. We can prove this equality using properties of the eta function: it satisfies the transformation formula
\begin{equation*}
    \eta(\gamma z) = \chi(\gamma) (-\i(cz+d))^{1/2} \eta(z)
\end{equation*}
for all $\gamma = \begin{psmallmatrix}
    a & b \\ c & d
\end{psmallmatrix} \in \SL_2(\Z)$ (with $c>0$) and $z \in \HH_2$, where 
\begin{equation*}
    \chi(\gamma) = \begin{cases}
        \exp\left( \pi \i \left( \tfrac{a+d}{12c} + s(-d,c) \right) \right) & c \ne 0, \\ 
        \exp\left( \tfrac{b\pi \i }{12} \right) & c = 0, d=1,
    \end{cases}
\end{equation*}
is a character of $\SL_2(\Z)$ with order 24, the function $s(\cdot, \cdot)$ is the Dedekind sum 
\begin{equation*}
    s(h,k) = \sum_{r=1}^{k-1} \frac{r}{k}\left( \frac{hr}{k} - \left[ \frac{hr}{k} \right] -\frac{1}{2}\right),
\end{equation*}
and $\left[ n \right]$ denotes the largest integer less than $n$. For a proof, see Theorem 3.4 of \cite{ApostolBook}. Using this formula on each factor on the left of \cref{Equation42EtaRelation}, we see that 
\begin{equation*}
    \eta(w) = \eta\left( \begin{pmatrix}
        1 & -11 \\ 2 & -21
    \end{pmatrix}\cdot (21w)  \right) = \exp(-5 \pi \i/6) (21)^{\frac{1}{4}} \eta(21w),
\end{equation*}

\begin{equation*}
    \eta(6w) =  \eta\left( \begin{pmatrix}
        3 & -22 \\ 1 & -7
    \end{pmatrix}\cdot (14w)  \right) = \exp(-\pi \i /3) \left( \frac{7}{3}\right)^{\frac{1}{4}} \eta(14w).
\end{equation*}

\begin{equation*}
    \eta(7w) = \eta\left( \begin{pmatrix}
        7 & -11 \\ 2 & -3
    \end{pmatrix}\cdot (3w)  \right) =  \exp(\pi \i /6)  \left( \frac{7}{3}\right)^{-\frac{1}{4}} \eta(3w),
\end{equation*}

\begin{equation*}
    \eta(42w) = \eta\left( \begin{pmatrix}
        21 & -22 \\ 1 & -1
    \end{pmatrix}\cdot (2w)  \right) = \exp(5\pi \i/3) (21)^{-\frac{1}{4}} \eta(2w), 
\end{equation*}

and so 
\begin{equation*}
    \eta(w) \eta(6w) \eta(7w) \eta(42w) = \exp(2 \pi \i /3) \eta(2w)\eta(3w)\eta(14w)\eta(21w).
\end{equation*}
Cubing both sides, we obtain \cref{Equation42EtaRelation}.

\section{Table of expressions} \label{SectionTables}

Forms are referred to with their LMFDB labels. Rows where we couldn't find an expression are left in the table, with a bound given on the smallest $d$ needed to raise the level of the form until it becomes expressible in terms of eta. For example, the form \LMRefDB{https://www.lmfdb.org/ModularForm/GL2/Q/holomorphic/53/2/a/a/}{53.2.a.a} has no (holomorphic) eta expression for all levels $\Gamma_0(53d)$, where $d=1,2,\ldots,17$, which is indicated by $18 \leq d$ in the ``Quotients'' column of the table. 

When an expression is found but has an unwieldy number of terms, we indicate the level of the form in the ``Quotients'' column, and in the ``Minimal'' column we either confirm minimality or bound the size $n$ of a minimal expression. For example, both forms of level $37$ are only expressible in terms of holomorphic eta quotients at level $37\cdot 8$, and have expressions of length 21 and 16 respectively. We have also confirmed that no expression exists for either form up to length $n=6$, and the ``Minimal'' column indicates this with $7 \leq n \leq 21$ and $7 \leq n \leq 16$ respectively. In the case of expressions not confirmed to be minimal, we always have the trivial bound $2 \leq d$, since the list of Martin and Ono \cite{MartinOno} is complete. 

Finally, in the ``Source'' column, we note the source of an eta expression. Some expressions are already found in the literature, and we include those that are minimal. Those with no source are original. 

An eta expression is primarily useful to us as a way of writing down a modular form quickly and easily. In cases where the expression is long (such as in the case of the two forms of level $37$), we omit the expression, as its length defeats this purpose. The expressions we found in these cases can be accessed at
\begin{center}
    \url{https://github.com/lewismcombes/EtaExpressions}
\end{center}

\newpage
\begin{table}[H]
\centering
\renewcommand*{\arraystretch}{1.1}
\begin{tabular}{|c|c|c|c|c|c|}\hline 
$N$  & LMFDB & Quotients & Coefficients & Minimal & Source \\ \hline 

11
&
\LMRefDB{https://www.lmfdb.org/ModularForm/GL2/Q/holomorphic/11/2/a/a/}{11.2.a.a}
&
$\eta_{11}[2,2]$
&
1
&
\checkmark
&
\cite{MartinOno}
\\ \hline
14
&
\LMRefDB{https://www.lmfdb.org/ModularForm/GL2/Q/holomorphic/14/2/a/a/}{14.2.a.a}
&
$\eta_{14}[1,1,1,1]$
&
1
& 
\checkmark
&
\cite{MartinOno}
\\ \hline
15
&
\LMRefDB{https://www.lmfdb.org/ModularForm/GL2/Q/holomorphic/15/2/a/a/}{15.2.a.a}
&
$\eta_{15}[1,1,1,1]$
&
1
&
\checkmark
&
\cite{MartinOno}
\\ \hline
17
&
\LMRefDB{https://www.lmfdb.org/ModularForm/GL2/Q/holomorphic/17/2/a/a/}{17.2.a.a}
&
$
\begin{matrix}
\eta_{68}[1,-1,2,-1,5,-2]\\ 
\eta_{68}[-1,5,-2,1,-1,2]
\end{matrix}
$
&
$
\begin{matrix} 
1 \\ -1 
\end{matrix}
$
&
\checkmark 
&
\cite{AyginThesis} 
\\ \hline
19
&
\LMRefDB{https://www.lmfdb.org/ModularForm/GL2/Q/holomorphic/19/2/a/a/}{19.2.a.a}
&
$
\begin{matrix}
\eta_{76}[ 4, -2, 0, 4, -2, 0 ]\\
\eta_{76}[ 3, -1, 0, -1, 3, 0 ]\\ 
\eta_{76}[ 1, 0, 1, -3, 8, -3 ]\\
\eta_{76}[ 2, -2, 2, 2, -2, 2 ]\\
\eta_{76}[ 0, 4, -2, 0, 4, -2 ]\\
\eta_{76}[ 0, 3, -1, 0, -1, 3 ]\\
\end{matrix}
$
&
$
\begin{matrix} 
1/4 \\ -2 \\ 2 \\ 1 \\ -1/4 \\ 4
\end{matrix}
$
&
\checkmark
&
\\ \hline
20
&
\LMRefDB{https://www.lmfdb.org/ModularForm/GL2/Q/holomorphic/20/2/a/a/}{20.2.a.a}
&
$\eta_{20}[0,2,0,0,2,0]$
&
1
&
\checkmark
&
\cite{MartinOno}
\\ \hline
21
&
\LMRefDB{https://www.lmfdb.org/ModularForm/GL2/Q/holomorphic/21/2/a/a/}{21.2.a.a}
&
$
\begin{matrix}
\eta_{42}[1,0,-1,2,0,1,2,-1] \\
\eta_{42}[0,1,2,-1,1,0,-1,2] \\
\eta_{42}[2,-1,-1,2,2,-1,-1,2]
\end{matrix}
$
&
$
\begin{matrix}
1 \\ 1 \\ 1
\end{matrix}
$
&
\checkmark
& 
\\ \hline
24 
&
\LMRefDB{https://www.lmfdb.org/ModularForm/GL2/Q/holomorphic/24/2/a/a/}{24.2.a.a}
&
$\eta_{24}[0,1,0,1,1,0,1,0]$
&
1 
&
\checkmark 
&
\cite{MartinOno}
\\ \hline 
\multirow{2}{*}{26} 
& 
\LMRefDB{https://www.lmfdb.org/ModularForm/GL2/Q/holomorphic/26/2/a/a/}{26.2.a.a}
& 
$
\begin{matrix} 
\eta_{52}[ -4, 10, -4, -4, 10, -4 ] \\ 
\eta_{52}[ -2, 7, -3, 2, 1, -1 ] \\
\eta_{52}[ -3, 7, -2, -1, 1, 2 ] \\
\eta_{52}[ 3, -2, 1, 1, -2, 3 ] \\
\eta_{52}[ 2, 1, -1, -2, 7, -3 ] \\
\eta_{52}[ 0, 0, 0, 0, -4, 8 ] \\
\eta_{52}[ 1, -2, 3, 3, -2, 1 ] \\ 
\eta_{52}[ 0, -4, 8, 0, 0, 0 ]\\
\eta_{52}[ 8, -4, 0, 0, 0, 0 ]
\end{matrix} 
$
& 
$
\begin{matrix} 
-7/12 \\ 3/4 \\ -4 \\ -5/3 \\ -1/4 \\ -26/3 \\ 1/3 \\ 2 \\ 1/12
\end{matrix}  
$
&
\checkmark
&
\\ 
\cline{2-6}
& 
\LMRefDB{https://www.lmfdb.org/ModularForm/GL2/Q/holomorphic/26/2/a/b/}{26.2.a.b}
& 
$
\begin{matrix}
\eta_{52}[ -2, 7, -3, 2, 1, -1 ] \\
\eta_{52}[ 3, -2, 1, 1, -2, 3 ] \\ 
\eta_{52}[ 2, 1, -1, -2, 7, -3 ] \\ 
\eta_{52}[ 1, -2, 3, 3, -2, 1 ] 
\end{matrix}
$
& 
$
\begin{matrix}
1/4 \\ 1 \\ -1/4 \\ 1
\end{matrix} 
$
&
\checkmark
&
\\
\hline

\end{tabular}
\end{table}


\newpage
\begin{table}[H]
\centering
\renewcommand*{\arraystretch}{1.1}
\begin{tabular}{|c|c|c|c|c|c|}\hline 
$N$  & LMFDB & Quotients & Coefficients & Minimal & Source \\ \hline 
27
&
\LMRefDB{https://www.lmfdb.org/ModularForm/GL2/Q/holomorphic/27/2/a/a/}{27.2.a.a}
&
$\eta_{27}[0,2,2,0]$
&
1
&
\checkmark
&
\cite{MartinOno}
\\ \hline
30
&
\LMRefDB{https://www.lmfdb.org/ModularForm/GL2/Q/holomorphic/30/2/a/a/}{30.2.a.a}
&
$
\begin{matrix}
\eta_{30}[ 0, 0, 1, 1, 1, 1, 0, 0 ] \\
\eta_{30}[ 1, 1, 0, 0, 0, 0, 1, 1 ]
\end{matrix}
$
&
$
\begin{matrix}
1 \\ -1
\end{matrix}
$
&
\checkmark
&
\\ \hline
32
&
\LMRefDB{https://www.lmfdb.org/ModularForm/GL2/Q/holomorphic/32/2/a/a/}{32.2.a.a}
&
$\eta_{32}[0,0,2,2,0,0]$
&
1
&
\checkmark
&
\cite{MartinOno}
\\ \hline
33
&
\LMRefDB{https://www.lmfdb.org/ModularForm/GL2/Q/holomorphic/33/2/a/a/}{33.2.a.a}
&
$
\begin{matrix}
\eta_{33}[2,0,2,0] \\
\eta_{33}[0,2,0,2] \\
\eta_{33}[1,1,1,1]
\end{matrix}
$
&
$
\begin{matrix}
1 \\ 3 \\ 3
\end{matrix}
$
&
\checkmark
&
\cite{PRY}
\\ \hline
34
&
\LMRefDB{https://www.lmfdb.org/ModularForm/GL2/Q/holomorphic/34/2/a/a/}{34.2.a.a}
&
$
\begin{matrix}
\eta_{68}[1,2,-1,-1,2,1] \\
\eta_{68}[-2,5,-1,2,-1,1] \\
\eta_{68}[-1,2,1,1,2,-1] \\
\eta_{68}[2,-1,1,-2,5,-1]
\end{matrix}
$
&
$
\begin{matrix}
1 \\ -1 \\ 1 \\ -1
\end{matrix}
$
&
\checkmark
&
\\ \hline
35
&
\LMRefDB{https://www.lmfdb.org/ModularForm/GL2/Q/holomorphic/35/2/a/a/}{35.2.a.a}
&
$
\begin{matrix}
\eta_{35}[2,0,0,2] \\
\eta_{35}[0,2,2,0] 
\end{matrix}
$
&
$
\begin{matrix}
1 \\ 1
\end{matrix}
$
&
\checkmark
&
\cite{AllenEtAl2020}
\\ \hline
36
&
\LMRefDB{https://www.lmfdb.org/ModularForm/GL2/Q/holomorphic/36/2/a/a/}{36.2.a.a}
&
$\eta_{36}[0,0,0,0,4,0,0,0,0]$
&
1
&
\checkmark
&
\cite{MartinOno}
\\ \hline
\multirow{2}{*}{37} 
& 
\LMRefDB{https://www.lmfdb.org/ModularForm/GL2/Q/holomorphic/37/2/a/a/}{37.2.a.a}
& 
$d = 8 = 2^3$
&
-
&
$7 \leq n \leq 21$
&
\\ 
\cline{2-6}
& 
\LMRefDB{https://www.lmfdb.org/ModularForm/GL2/Q/holomorphic/37/2/a/b/}{37.2.a.b}
& 
$d = 8 = 2^3$
&
-
&
$7\leq n \leq 16$
&
\\
\hline
\multirow{2}{*}{38} 
& 
\LMRefDB{https://www.lmfdb.org/ModularForm/GL2/Q/holomorphic/38/2/a/a/}{38.2.a.a}
& 
$
\begin{matrix}
\eta_{76}[ -1, 3, 0, 3, -1, 0 ] \\
\eta_{76}[ -2, 4, 0, -2, 4, 0 ] \\
\eta_{76}[ 4, -2, 0, 4, -2, 0 ] \\
\eta_{76}[ 3, -1, 0, -1, 3, 0 ] \\
\eta_{76}[ -4, 10, -4, -4, 10, -4 ] \\
\eta_{76}[ 0, -2, 4, 0, -2, 4 ] 
\end{matrix}
$
& 
$
\begin{matrix}
-1 \\ -2 \\ -1/4 \\ 1 \\ 1/4 \\ 4
\end{matrix}
$
&
\checkmark
&
\\ 
\cline{2-6}
& 
\LMRefDB{https://www.lmfdb.org/ModularForm/GL2/Q/holomorphic/38/2/a/b/}{38.2.a.b}
& 
$
\begin{matrix}
\eta_{38}[-1,3,3,-1] \\
\eta_{38}[3,-1,-1,3]
\end{matrix}$
& 
$
\begin{matrix}
1 \\ 1
\end{matrix}
$
&
\checkmark
&
\cite{AyginThesis}
\\
\hline
\end{tabular}
\end{table}


\newpage
\begin{table}[H]
\centering
\renewcommand*{\arraystretch}{1.1}
\begin{tabular}{|c|c|c|c|c|c|}\hline 
$N$  & LMFDB & Quotients & Coefficients & Minimal & Source \\ \hline 
39
&
\LMRefDB{https://www.lmfdb.org/ModularForm/GL2/Q/holomorphic/39/2/a/a/}{39.2.a.a}
&
$
\begin{matrix}
\eta_{78}[ 0, 0, 1, 1, 1, 1, 0, 0 ] \\
\eta_{78}[ 0, 0, 4, -2, 4, -2, 0, 0 ] \\
\eta_{78}[ 4, -2, 0, 0, 0, 0, 4, -2 ] \\
\eta_{78}[ 1, 1, 0, 0, 0, 0, 1, 1 ]
\end{matrix}
$
&
$
\begin{matrix}
2 \\ 1/4 \\ -1/4 \\ 2
\end{matrix}
$
&
\checkmark
&
\\ \hline
40
&
\LMRefDB{https://www.lmfdb.org/ModularForm/GL2/Q/holomorphic/40/2/a/a/}{40.2.a.a}
&
$
\begin{matrix}
\eta_{40}[-1,2,2,1,-1,0,0,1] \\ 
\eta_{40}[1,0,0,-1,1,2,2,-1]
\end{matrix}
$
&
$
\begin{matrix}
1 \\ 1
\end{matrix}
$
&
\checkmark
&
\cite{PRY}
\\ \hline
42
&
\LMRefDB{https://www.lmfdb.org/ModularForm/GL2/Q/holomorphic/42/2/a/a/}{42.2.a.a}
&
$
\begin{matrix}
\eta_{42}[ -1, 2, 2, -1, -1, 2, 2, -1 ] \\ 
\eta_{42}[ 2, -1, -1, 2, 2, -1, -1, 2 ]
\end{matrix}
$
&
$
\begin{matrix}
1 \\ -1
\end{matrix}
$
&
\checkmark
&
\\ \hline
43
&
\LMRefDB{https://www.lmfdb.org/ModularForm/GL2/Q/holomorphic/43/2/a/a/}{43.2.a.a}
&
$d= 8 = 2^3$
&
-
&
$6 \leq n \leq 20$
&

\\ \hline
44
&
\LMRefDB{https://www.lmfdb.org/ModularForm/GL2/Q/holomorphic/44/2/a/a/}{44.2.a.a}
&
$
\begin{matrix}
\eta_{44}[1,-1,0,-3,7,0] \\ 
\eta_{44}[0,-1,1,0,7,-3] \\
\eta_{44}[-1,2,-1,3,-2,3] 
\end{matrix}
$
&
$
\begin{matrix}
-2 \\ 1 \\ -2 
\end{matrix}
$
&
\checkmark
&
\\ \hline
45
&
\LMRefDB{https://www.lmfdb.org/ModularForm/GL2/Q/holomorphic/45/2/a/a/}{45.2.a.a}
&
$
\begin{matrix}
\eta_{45}[1,0,-2,1,6,-2] \\ 
\eta_{45}[0,0,-3,0,10,-3] \\
\eta_{45}[-1,2,2,-1,0,2] 
\end{matrix}
$
&
$
\begin{matrix}
-1 \\ 1 \\ -1 
\end{matrix}
$
&
\checkmark
&
\\ \hline
46
&
\LMRefDB{https://www.lmfdb.org/ModularForm/GL2/Q/holomorphic/46/2/a/a/}{46.2.a.a}
&
$
\begin{matrix}
\eta_{46}[3,-1,3,-1] \\ 
\eta_{46}[2,0,2,0] \\
\eta_{46}[1,1,1,1] \\ 
\eta_{46}[0,2,0,2] \\
\eta_{46}[-1,3,-1,3]
\end{matrix}
$
&
$
\begin{matrix}
1 \\ 2 \\ 3 \\ 4 \\ 4
\end{matrix}
$
&
\checkmark
&

\\ \hline
48
&
\LMRefDB{https://www.lmfdb.org/ModularForm/GL2/Q/holomorphic/48/2/a/a/}{48.2.a.a}
&
$\eta_{48}[0,-1,0,4,-1,-1,4,0,-1,0]$
&
1
&
\checkmark
&
\cite{MartinOno}
\\ \hline
49
&
\LMRefDB{https://www.lmfdb.org/ModularForm/GL2/Q/holomorphic/49/2/a/a/}{49.2.a.a}
&
$d = 4 = 2^2$
&
-
&
$4 \leq n \leq 8$
&
\\ \hline
\end{tabular}
\end{table}


\newpage
\begin{table}[H]
\centering
\renewcommand*{\arraystretch}{1.1}
\begin{tabular}{|c|c|c|c|c|c|}\hline 
$N$  & LMFDB & Quotients & Coefficients & Minimal & Source \\ \hline 

\multirow{2}{*}{50} 
& 
\LMRefDB{https://www.lmfdb.org/ModularForm/GL2/Q/holomorphic/50/2/a/a/}{50.2.a.a}
& 
$
\begin{matrix}
\eta_{100}[ -1, 3, -1, 0, 2, 0, 1, -1, 1 ]\\
\eta_{100}[ 1, -1, 1, 0, 2, 0, -1, 3, -1 ] \\
\eta_{100}[ 0, 0, 1, 0, 2, 0, 0, 2, -1 ] \\
\eta_{100}[ 0, 2, -1, 0, 2, 0, 0, 0, 1 ]
\end{matrix}$
& 
$
\begin{matrix}
1 \\ -1 \\ 1 \\ -1
\end{matrix}
$
&
\checkmark
&
\\ 
\cline{2-6}
& 
\LMRefDB{https://www.lmfdb.org/ModularForm/GL2/Q/holomorphic/50/2/a/b/}{50.2.a.b}
& 
$
\begin{matrix}
\eta_{100}[ -1, 2, -1, 2, 0, 2, -1, 2, -1 ] \\
\eta_{100}[ 0, -1, 1, 0, 6, -2, 0, -1, 1 ]
\end{matrix}$
&
$
\begin{matrix}
1 \\ -1
\end{matrix}
$
&
\checkmark 
&
\\
\hline
51
&
\LMRefDB{https://www.lmfdb.org/ModularForm/GL2/Q/holomorphic/51/2/a/a/}{51.2.a.a}
&
$d = 2$
&
-
&
$6 \leq n \leq 9$
&
\\ \hline

52
&
\LMRefDB{https://www.lmfdb.org/ModularForm/GL2/Q/holomorphic/52/2/a/a/}{52.2.a.a}
&
$ 
\begin{matrix}
\eta_{52}[-4,10,-4,-4,10,-4] \\
\eta_{52}[4,-2,0,4,-2,0]\\
\eta_{52}[0,-2,4,0,-2,4]
\end{matrix}$
&
$
\begin{matrix}
1/8 \\ -1/8 \\ -2
\end{matrix}
$
&
\checkmark
&
\cite{AyginThesis}
\\ \hline
53
&
\LMRefDB{https://www.lmfdb.org/ModularForm/GL2/Q/holomorphic/53/2/a/a/}{53.2.a.a}
&
$ 18 \leq d $
&
-
&

&
\\ \hline
\multirow{2}{*}{54} 
& 
\LMRefDB{https://www.lmfdb.org/ModularForm/GL2/Q/holomorphic/54/2/a/a/}{54.2.a.a}
& 
$
\begin{matrix}
\eta_{54}[ 0, 0, -1, 3, 3, -1, 0, 0 ]\\
\eta_{54}[ 0, 0, 3, -1, -1, 3, 0, 0 ]
\end{matrix}$
& 
$
\begin{matrix}
1 \\ -1
\end{matrix}
$
&
\checkmark 
&

\\ 
\cline{2-6}
& 
\LMRefDB{https://www.lmfdb.org/ModularForm/GL2/Q/holomorphic/54/2/a/b/}{54.2.a.b}
& 
$
\begin{matrix}
\eta_{54}[ 0, 0, -1, 3, 3, -1, 0, 0 ]\\
\eta_{54}[ 0, 0, 3, -1, -1, 3, 0, 0 ]
\end{matrix}$
& 
$
\begin{matrix}
1 \\ 1
\end{matrix}
$
&
\checkmark
&
\\
\hline
55
&
\LMRefDB{https://www.lmfdb.org/ModularForm/GL2/Q/holomorphic/55/2/a/a/}{55.2.a.a}
&
$
\begin{matrix}
\eta_{110}[ 0, 0, 1, 1, 1, 1, 0, 0 ] \\ 
\eta_{110}[ 1, -1, -1, 3, -1, 3, 1, -1 ] \\ 
\eta_{110}[ -1, 3, 1, -1, 1, -1, -1, 3 ] \\ 
\eta_{110}[0,0,4,-2,1,-1,-1,3] \\
\end{matrix}$
&
$
\begin{matrix}
2 \\ 1 \\ -1 \\  2
\end{matrix}$
&
\checkmark
&
\\ \hline
\multirow{2}{*}{56} 
& 
\LMRefDB{https://www.lmfdb.org/ModularForm/GL2/Q/holomorphic/56/2/a/a/}{56.2.a.a}
& 
$
\begin{matrix}
\eta_{56}[ 0, -1, 3, 0, 0, 3, -1, 0 ] \\
\eta_{56}[ 0, 3, -1, 0, 0, -1, 3, 0 ] 
\end{matrix}$
& 
$
\begin{matrix}
1 \\ -1
\end{matrix}
$
&
\checkmark
&

\\ 
\cline{2-6}
& 
\LMRefDB{https://www.lmfdb.org/ModularForm/GL2/Q/holomorphic/56/2/a/b/}{56.2.a.b}
& 
$
\begin{matrix}
\eta_{56}[ 0, -1, 3, 0, 0, 3, -1, 0 ] \\
\eta_{56}[ 0, 3, -1, 0, 0, -1, 3, 0 ]
\end{matrix}$
& 
$
\begin{matrix}
1 \\ 1
\end{matrix}
$
&
\checkmark
&

\\
\hline
\multirow{2}{*}{57} 
& 
\LMRefDB{https://www.lmfdb.org/ModularForm/GL2/Q/holomorphic/57/2/a/a/}{57.2.a.a}
& 
$d = 2$
& 
-
&
$5 \leq n \leq 7$
&
\\ 
\cline{2-6}
& 
\LMRefDB{https://www.lmfdb.org/ModularForm/GL2/Q/holomorphic/57/2/a/b/}{57.2.a.b}
& 
$d = 2$
& 
-
&
$5 \leq n \leq 8$
&
\\
\cline{2-6}
& 
\LMRefDB{https://www.lmfdb.org/ModularForm/GL2/Q/holomorphic/57/2/a/c/}{57.2.a.c}
& 
$d = 2$
& 
-
&
$6 \leq n \leq 8$
&
\\
\hline
\end{tabular}
\end{table}


\newpage
\begin{table}[H]
\centering
\renewcommand*{\arraystretch}{1.1}
\begin{tabular}{|c|c|c|c|c|c|}\hline 
$N$  & LMFDB & Quotients & Coefficients & Minimal & Source \\ \hline 
\multirow{2}{*}{58} 
& 
\LMRefDB{https://www.lmfdb.org/ModularForm/GL2/Q/holomorphic/58/2/a/a/}{58.2.a.a}
& 
$
\begin{matrix}
\eta_{116}[ -4, 10, -4, -4, 10, -4 ] \\
\eta_{116}[ -2, 7, -3, 2, 1, -1 ] \\
\eta_{116}[ 3, -2, 1, 1, -2, 3 ] \\
\eta_{116}[ 2, 1, -1, -2, 7, -3 ] \\
\eta_{116}[ 0, 0, 0, 8, -4, 0 ] \\
\eta_{116}[ 0, 0, 0, 0, -4, 8 ] \\
\eta_{116}[ 1, -2, 3, 3, -2, 1 ] \\
\eta_{116}[ 0, -2, 4, 0, -2, 4 ]  \\
\eta_{116}[ 0, -4, 8, 0, 0, 0 ] \\
\eta_{116}[ 8, -4, 0, 0, 0, 0 ]
\end{matrix}$
& 
$
\begin{matrix}
-5/16 \\  1/2  \\   2 \\ -1/2\\ 29/96 \\ 29/3  \\   2  \\  -5 \\  1/3 \\ 1/96
\end{matrix}
$
&
\checkmark
&

\\ 
\cline{2-6}
& 
\LMRefDB{https://www.lmfdb.org/ModularForm/GL2/Q/holomorphic/58/2/a/b/}{58.2.a.b}
& 
$
\begin{matrix}
\eta_{116} [ -4, 10, -4, -4, 10, -4 ] \\
\eta_{116}[ 0, 0, 0, 8, -4, 0 ] \\
\eta_{116}[ 0, 0, 0, 0, -4, 8 ] \\
\eta_{116}[ 0, -2, 4, 0, -2, 4 ] \\
\eta_{116}[ 0, -4, 8, 0, 0, 0 ] \\
\eta_{116}[ 8, -4, 0, 0, 0, 0 ] 
\end{matrix}$
& 
$
\begin{matrix}
5/16 \\ -29/96 \\  -29/3  \\    5 \\  -1/3 \\ -1/96
\end{matrix}
$
&
\checkmark
&
\\
\hline

61
&
\LMRefDB{https://www.lmfdb.org/ModularForm/GL2/Q/holomorphic/61/2/a/a/}{61.2.a.a}
&
$d = 12 = 2^2 \cdot 3$
&
-
&
$2 \leq n \leq 101$
&
\\ \hline
62
&
\LMRefDB{https://www.lmfdb.org/ModularForm/GL2/Q/holomorphic/62/2/a/a/}{62.2.a.a}
&
$d = 3$
&
-
&
$7 \leq  n \leq 8$
&
\\ \hline
63
&
\LMRefDB{https://www.lmfdb.org/ModularForm/GL2/Q/holomorphic/63/2/a/a/}{63.2.a.a}
&
$
\begin{matrix}
\eta_{63}[ 0, 0, 3, 0, -2, 3 ] \\
\eta_{63}[ 3, -2, 0, 3, 0, 0 ] \\
\eta_{63}[ 0, -1, 3, 3, -1, 0 ] \\
\eta_{63}[ 3, -1, 0, 0, -1, 3 ]
\end{matrix}
$
&
$
\begin{matrix}
-7/3 \\ -1/3 \\  4/3 \\  4/3
\end{matrix}$
&
\checkmark
&
\\ \hline
64
&
\LMRefDB{https://www.lmfdb.org/ModularForm/GL2/Q/holomorphic/64/2/a/a/}{64.2.a.a}
&
$\eta_{64}[0,0,-2,8,-2,0,0]$
&
1
&
\checkmark
&
\cite{MartinOno}
\\ \hline

\end{tabular}
\end{table}


\newpage
\begin{table}[H]
\centering
\renewcommand*{\arraystretch}{1.1}
\begin{tabular}{|c|c|c|c|c|c|}\hline 
$N$  & LMFDB & Quotients & Coefficients & Minimal & Source \\ \hline 
65
&
\LMRefDB{https://www.lmfdb.org/ModularForm/GL2/Q/holomorphic/65/2/a/a/}{65.2.a.a}
&
$d = 4 = 2^2$
&
-
&
$4 \leq n \leq 6$
&
\\ \hline
\multirow{2}{*}{66} 
& 
\LMRefDB{https://www.lmfdb.org/ModularForm/GL2/Q/holomorphic/66/2/a/a/}{66.2.a.a}
& 
$
\begin{matrix}
\eta_{66}[-3,6,1,-2,-2,1,6,-3] \\
\eta_{66}[-2,1,6,-3,-3,6,1,-2]
\end{matrix}$
& 
$
\begin{matrix}
1 \\ -1
\end{matrix}
$
&
\checkmark
&
\\ 
\cline{2-6}
& 
\LMRefDB{https://www.lmfdb.org/ModularForm/GL2/Q/holomorphic/66/2/a/b/}{66.2.a.b}
& 
$
\begin{matrix}
\eta_{66}[1,2,0,-1,-1,0,2,1] \\ 
\eta_{66}[5,-2,-1,0,0,-1,-2,5] \\ 
\eta_{66}[-1,0,2,1,1,2,0,-1] \\
\eta_{66}[0,-1,-2,5,5,-2,-1,0]
\end{matrix}$
& 
$
\begin{matrix}
1 \\ 1 \\ 1 \\ -1
\end{matrix}
$
&
\checkmark
&
\\
\cline{2-6}
& 
\LMRefDB{https://www.lmfdb.org/ModularForm/GL2/Q/holomorphic/66/2/a/c/}{66.2.a.c}
& 
$
\begin{matrix}
\eta_{66}[-1,3,1,-1,-1,3,1,-1] \\
\eta_{66}[-1,2,1,0,-1,2,1,0] \\ 
\eta_{66}[1,0,-1,2,1,0,-1,2] \\
\eta_{66}[1,-1,-1,3,1,-1,-1,3]
\end{matrix}$
& 
$
\begin{matrix}
1 \\ 2 \\ -2 \\ -3
\end{matrix}
$
&
\checkmark
&
\\
\hline
67
&
\LMRefDB{https://www.lmfdb.org/ModularForm/GL2/Q/holomorphic/67/2/a/a/}{67.2.a.a}
&
$d = 12 = 2^2 \cdot 3$
&
-
&
$3 \leq n \leq 41$
&
\\ \hline
69
&
\LMRefDB{https://www.lmfdb.org/ModularForm/GL2/Q/holomorphic/69/2/a/a/}{69.2.a.a}
&
$d = 4 = 2^2$
&
-
&
$ 3 \leq n \leq 9$
&
\\ \hline
70
&
\LMRefDB{https://www.lmfdb.org/ModularForm/GL2/Q/holomorphic/70/2/a/a/}{70.2.a.a}
&
$
\begin{matrix}
\eta_{70}[-1,2,2,-1,-1,2,2,-1] \\
\eta_{70}[2,-1,-1,2,2,-1,-1,2]
\end{matrix}$
&
$
\begin{matrix}
1 \\ -1
\end{matrix}
$
&
\checkmark
&
\cite{PRY}
\\ \hline
72
&
\LMRefDB{https://www.lmfdb.org/ModularForm/GL2/Q/holomorphic/72/2/a/a/}{72.2.a.a}
&
$
\begin{matrix}
\eta_{72}[ 0, 0, 0, -2, 0, 4, 4, 0, -2, 0, 0, 0 ] \\ 
\eta_{72}[ 4, -2, 0, 0, 0, 0, 0, 0, 0, 0, -2, 4 ]
\end{matrix}$
&
$
\begin{matrix}
1 \\ -1 
\end{matrix}
$
&
\checkmark 
&

\\ \hline
73
&
\LMRefDB{https://www.lmfdb.org/ModularForm/GL2/Q/holomorphic/73/2/a/a/}{73.2.a.a}
&
$16 \leq d$
&
-
&

&
\\ \hline
\end{tabular}
\end{table}


\newpage
\begin{table}[H]
\centering
\renewcommand*{\arraystretch}{1.1}
\begin{tabular}{|c|c|c|c|c|c|}\hline 
$N$  & LMFDB & Quotients & Coefficients & Minimal & Source \\ \hline 
\multirow{2}{*}{75} 
& 
\LMRefDB{https://www.lmfdb.org/ModularForm/GL2/Q/holomorphic/75/2/a/a/}{75.2.a.a}
& 
$
\begin{matrix}
\eta_{75}[ 0, 0, -2, 4, 4, -2 ] \\
\eta_{75}[ 1, 1, 1, 1, 0, 0 ] \\
\eta_{75}[ 1, 0, 1, 1, 0, 1 ] \\
\eta_{75}[ 0, 1, 1, 1, 1, 0 ] \\
\eta_{75}[ -1, 2, 1, 1, 2, -1 ] \\
\eta_{75}[ 0, 0, 1, 1, 1, 1 ] \\
\eta_{75}[ -1, 3, -1, 3, 0, 0 ] \\
\eta_{75}[ 0, 0, -1, 3, -1, 3 ]
\end{matrix}$
& 
$
\begin{matrix}
 -4  \\  -3 \\ -1\\ -7 \\ 4 \\   -5  \\ -6 \\ 6
\end{matrix}
$
&
\checkmark
&
\\ 
\cline{2-6}
& 
\LMRefDB{https://www.lmfdb.org/ModularForm/GL2/Q/holomorphic/75/2/a/b/}{75.2.a.b}
& 
$
\begin{matrix}
\eta_{75}[ 2, -1, 1, 1, -1, 2 ] \\
\eta_{75}[ 0, 0, 4, -2, -2, 4 ] \\
\eta_{75}[ 3, -1, 3, -1, 0, 0 ] \\
\eta_{75}[ -1, 3, -1, 3, 0, 0 ] \\
\eta_{75}[ 0, 0, -2, 4, 4, -2 ] \\
\eta_{75} [ -1, 2, 1, 1, 2, -1 ] \\
\eta_{75}[ 0, 0, 3, -1, 3, -1 ] \\
\eta_{75}[ 0, 0, -1, 3, -1, 3 ]
\end{matrix}$
& 
$
\begin{matrix}
 -2 \\  -2 \\ 1/3 \\  -3 \\  -2  \\  2\\ -1/3 \\   3
\end{matrix}
$
&
\checkmark
&
\\
\cline{2-6}
& 
\LMRefDB{https://www.lmfdb.org/ModularForm/GL2/Q/holomorphic/75/2/a/c/}{75.2.a.c}
& 
$
\begin{matrix}
\eta_{75}[0,1,1,1,1,0] \\ 
\eta_{75}[0,0,1,1,1,1] \\
\eta_{75}[1,0,1,1,0,1] \\
\eta_{75}[1,1,1,1,0,0]
\end{matrix}$
& 
$
\begin{matrix}
3 \\ 5 \\ 3 \\ 1
\end{matrix}
$
&
\checkmark
&
\\
\hline
76
& 
\LMRefDB{https://www.lmfdb.org/ModularForm/GL2/Q/holomorphic/76/2/a/a/}{76.2.a.a}
& 
$
\begin{matrix}
\eta_{152}[ -2, 3, 3, -2, -2, 3, 3, -2 ] \\
\eta_{152}[ 2, -3, 5, -2, 2, -3, 5, -2 ] \\
\eta_{152}[ 0, -1, 1, 2, 0, -1, 1, 2 ] \\
\eta_{152}[ -2, 5, -3, 2, -2, 5, -3, 2 ] \\
\eta_{152}[ 2, -1, -1, 2, 2, -1, -1, 2 ] \\
\eta_{152}[ 0, -3, 7, -2, 0, -3, 7, -2 ]
\end{matrix}$
& 
$
\begin{matrix}
1/4 \\ -1/4 \\ -4 \\ 1 \\ -1 \\ -1
\end{matrix}
$
&
$n=6$
&
\\ \hline
\end{tabular}
\end{table}


\newpage
\begin{table}[H]
\centering
\renewcommand*{\arraystretch}{1.1}
\begin{tabular}{|c|c|c|c|c|c|}\hline 
$N$  & LMFDB & Quotients & Coefficients & Minimal & Source \\ \hline 

\multirow{2}{*}{77} 
& 
\LMRefDB{https://www.lmfdb.org/ModularForm/GL2/Q/holomorphic/77/2/a/a/}{77.2.a.a}
& 
$d = 4 = 2^2$
& 
-
&
$4 \leq n \leq 8$
&
\\ 
\cline{2-6}
& 
\LMRefDB{https://www.lmfdb.org/ModularForm/GL2/Q/holomorphic/77/2/a/b/}{77.2.a.b}
& 
$
\begin{matrix}
\eta_{154}[2,-1,2,-1,-1,2,-1,2] \\
\eta_{154}[2,-1,-1, 0, 2,1,1,0]\\
\eta_{154}[1,0,0,-1,1,2,2,-1]\\
\eta_{154}[1,0,-1,1,2,0,-1,2]\\
\eta_{154}[-1,2,1,-1,0,2,1,0] \\
\eta_{154}[1,0,1,1,0,0,1,0] \\
\eta_{154}[0,1,1,2,0,-1,-1,2] \\
\eta_{154}[-1,2,2,1,-1,0,0,1] \\
\eta_{154}[-1,2,-1,2,2,-1,2,-1]
\end{matrix}$
& 
$
\begin{matrix}
1 \\ 1 \\ -1 \\ -2 \\ -2 \\ -4 \\ -1 \\ 1 \\ 1
\end{matrix}
$
&
\checkmark
&
\\
\cline{2-6}
& 
\LMRefDB{https://www.lmfdb.org/ModularForm/GL2/Q/holomorphic/77/2/a/c/}{77.2.a.c}
& 
$d = 4 = 2^2$
& 
-
&
$4 \leq n \leq 6$
&
\\
\hline
78 
&
\LMRefDB{https://www.lmfdb.org/ModularForm/GL2/Q/holomorphic/78/2/a/a/}{78.2.a.a}
& 
$
\begin{matrix}
\eta_{78}[-3,6,1,-2,-2,1,6,-3]\\
\eta_{78}[1,0,0,1,0,1,1,0] \\ 
\eta_{78}[0,1,1,0,1,0,0,1] \\
\eta_{78}[-2,1,6,-3,-3,6,1,-2]
\end{matrix}
$
&
$
\begin{matrix}
1 \\ -2 \\ -2 \\ -1
\end{matrix}
$
& 
\checkmark
&
 \\ \hline 
79
&
\LMRefDB{https://www.lmfdb.org/ModularForm/GL2/Q/holomorphic/79/2/a/a/}{79.2.a.a}
& 
$d = 12 = 2^2 \cdot 3$
&
-
& 
$3 \leq n \leq 155$
&
 \\ \hline 
\multirow{2}{*}{80} 
& 
\LMRefDB{https://www.lmfdb.org/ModularForm/GL2/Q/holomorphic/80/2/a/a/}{80.2.a.a}
& 
$
\begin{matrix}
\eta_{80}[ 0, 0, -2, 0, 4, 4, 0, -2, 0, 0 ] \\
\eta_{80}[ 0, 4, -2, 0, 0, 0, 0, -2, 4, 0 ] 
\end{matrix}$
& 
$
\begin{matrix}
1 \\ -1 
\end{matrix}
$
&
\checkmark
&

\\ 
\cline{2-6}
& 
\LMRefDB{https://www.lmfdb.org/ModularForm/GL2/Q/holomorphic/80/2/a/b/}{80.2.a.b}
& 
$\eta_{80}[0,-2,6,0,-2,-2,0,6,-2,0]$
& 
1
&
\checkmark 
&
\cite{MartinOno}
\\
\hline

82
&
\LMRefDB{https://www.lmfdb.org/ModularForm/GL2/Q/holomorphic/82/2/a/a/}{82.2.a.a}
& 
$d = 4 = 2^2$
&
-
& 
$6 \leq n \leq 14$
&
 \\ \hline 
83
&
\LMRefDB{https://www.lmfdb.org/ModularForm/GL2/Q/holomorphic/83/2/a/a/}{83.2.a.a}
& 
$12 \leq d$
&
-
& 

&
 \\ \hline 
\end{tabular}
\end{table}


\newpage
\begin{table}[H]
\centering
\renewcommand*{\arraystretch}{1.1}
\begin{tabular}{|c|c|c|c|c|c|}\hline 
$N$  & LMFDB & Quotients & Coefficients & Minimal & Source \\ \hline 

\multirow{2}{*}{84} 
& 
\LMRefDB{https://www.lmfdb.org/ModularForm/GL2/Q/holomorphic/84/2/a/a/}{84.2.a.a}
& 
$d = 1$
& 
-
&
$3\leq n \leq 7$
&
\\ 
\cline{2-6}
& 
\LMRefDB{https://www.lmfdb.org/ModularForm/GL2/Q/holomorphic/84/2/a/b/}{84.2.a.b}
& 
$\begin{matrix}
\eta_{84}[ 0, -1, 0, 2, -1, 2, 2, -1, 2, 0, -1, 0 ] \\
\eta_{84}[ 2, -1, 2, 0, -1, 0, 0, -1, 0, 2, -1, 2 ]
\end{matrix}$
& 
$
\begin{matrix}
1 \\ -1
\end{matrix}
$
&
\checkmark
&
\\
\hline
85
&
\LMRefDB{https://www.lmfdb.org/ModularForm/GL2/Q/holomorphic/85/2/a/a/}{85.2.a.a}
& 
$d = 4 = 2^2$
&
-
& 
$2 \leq n \leq 14$
&
 \\ \hline 
88
&
\LMRefDB{https://www.lmfdb.org/ModularForm/GL2/Q/holomorphic/88/2/a/a/}{88.2.a.a}
& 
$\begin{matrix}
\eta_{88}[0,3,-1,0,0,3,-1,0] \\ 
\eta_{88}[0,-1,3,0,0,-1,3,0]
\end{matrix}$

&
$
\begin{matrix}
1 \\ -4
\end{matrix}
$
& 
\checkmark
&
 \\ \hline 
\multirow{2}{*}{89} 
& 
\LMRefDB{https://www.lmfdb.org/ModularForm/GL2/Q/holomorphic/89/2/a/a/}{89.2.a.a}
& 
$18 \leq d$
& 
-
&

&
\\ 
\cline{2-6}
& 
\LMRefDB{https://www.lmfdb.org/ModularForm/GL2/Q/holomorphic/89/2/a/b/}{89.2.a.b}
& 
$d = 12 = 2^2 \cdot 3$
& 
-
&
$2 \leq n \leq 160$
&
\\
\hline
\multirow{2}{*}{90} 
& 
\LMRefDB{https://www.lmfdb.org/ModularForm/GL2/Q/holomorphic/90/2/a/a/}{90.2.a.a}
& 
$
\begin{matrix}
\eta_{90}[ 1, -1, -1, 0, 3, 1, 1, -1, -1, 1, 0, 1 ] \\ 
\eta_{90}[ 2, -1, -1, 0, 2, 0, 0, 1, 0, 0, 2, -1 ] \\
\eta_{90}[ 0, 0, 3, 2, -1, 0, -1, -3, 0, 3, 2, -1 ] \\
\eta_{90}[ -1, 0, 5, -1, -2, -2, 3, -1, 2, 0, 0, 1 ] 
\end{matrix}
$
& 
$
\begin{matrix}
2 \\ 1 \\ 1 \\ 2
\end{matrix}
$
&
$3 \leq n \leq 4$
&
\\ 
\cline{2-6}
& 
\LMRefDB{https://www.lmfdb.org/ModularForm/GL2/Q/holomorphic/90/2/a/b/}{90.2.a.b}
& 
$
\begin{matrix}
\eta_{90}[ 1, -1, -2, 0, 4, 2, 1, 3, -2, -3, -2, 3 ] \\ 
\eta_{90}[ 0, 0, 1, -1, 1, -1, 2, 2, 1, -2, 1, 0 ] \\
\eta_{90}[ 0, 0, 3, 0, -1, -1, 0, 0, 1, 0, 0, 2 ] \\
\eta_{90}[ 0, 0, -1, -1, 2, 2, 2, 2, -1, -1, 0, 0 ] 
\end{matrix}
$
& 
$
\begin{matrix}
-1 \\ 1 \\ 3 \\ 1
\end{matrix}
$
&
$3 \leq n \leq 4$
&
\\
\cline{2-6}
& 
\LMRefDB{https://www.lmfdb.org/ModularForm/GL2/Q/holomorphic/90/2/a/c/}{90.2.a.c}
& 
$
\begin{matrix}
\eta_{90}[ 0, 0, 2, 2, -1, -1, -1, -1, 2, 2, 0, 0 ] \\ 
\eta_{90}[ -2, 3, 2, 2, -2, -2, -2, -2, 3, 4, 2, -2 ] \\
\eta_{90}[ 0, 0, 1, -1, 0, 2, 2, 0, -1, 1, 0, 0 ] \\
\eta_{90}[ -1, 1, 0, 1, 2, -1, 0, 0, 1, 0, 1, 0 ] 
\end{matrix}
$
& 
$
\begin{matrix}
-1 \\ 1 \\ -1 \\ -1
\end{matrix}
$
&
$3 \leq n \leq 4$
&
\\
\hline
\multirow{2}{*}{91} 
& 
\LMRefDB{https://www.lmfdb.org/ModularForm/GL2/Q/holomorphic/91/2/a/a/}{91.2.a.a}
& 
$d = 4 = 2^2$
& 
-
&
$4 \leq n \leq 15$
&
\\ 
\cline{2-6}
& 
\LMRefDB{https://www.lmfdb.org/ModularForm/GL2/Q/holomorphic/91/2/a/b/}{91.2.a.b}
& 
$d = 4 = 2^2$
& 
-
&
$4 \leq n \leq 22$
&
\\
\hline

\end{tabular}
\end{table}


\newpage
\begin{table}[H]
\centering
\renewcommand*{\arraystretch}{1.1}
\begin{tabular}{|c|c|c|c|c|c|}\hline 
$N$  & LMFDB & Quotients & Coefficients & Minimal & Source \\ \hline 
\multirow{2}{*}{92} 
& 
\LMRefDB{https://www.lmfdb.org/ModularForm/GL2/Q/holomorphic/92/2/a/a/}{92.2.a.a}
& 
$
\begin{matrix}
\eta_{92}[ 2, -1, 1, 2, -1, 1 ] \\ 
\eta_{92}[ 1, -1, 2, 1, -1, 2 ] \\
\eta_{92}[ -1, 5, -2, -1, 5, -2 ] \\
\eta_{92}[ -2, 5, -1, -2, 5, -1 ] \\ 
\eta_{92}[ 0, 3, -1, 0, 3, -1 ] \\
\eta_{92}[ -1, 3, 0, -1, 3, 0 ]
\end{matrix}
$
& 
$
\begin{matrix}
1 \\ -2 \\ 1 \\ -2 \\ -1 \\ 2
\end{matrix}
$
&
\checkmark
&
\\ 
\cline{2-6}
& 
\LMRefDB{https://www.lmfdb.org/ModularForm/GL2/Q/holomorphic/92/2/a/b/}{92.2.a.b}
& 
$
\begin{matrix}
\eta_{276}[ 0, -2, 0, 3, 2, -1, 3, -2, -1, 0, 2, 0 ] \\ 
\eta_{276}[ 2, -1, 2, 0, -1, 0, 0, -1, 0, 2, -1, 2 ] \\
\eta_{276}[ 0, -1, 0, 2, -1, 2, 2, -1, 2, 0, -1, 0 ] \\
\eta_{276}[ 3, -2, -1, 0, 2, 0, 0, -2, 0, 3, 2, -1 ] 
\end{matrix}
$
& 
$
\begin{matrix}
-1 \\ -1 \\ 1 \\ 1
\end{matrix}
$
&
$3 \leq n \leq 4$
&
\\
\hline
94
& 
\LMRefDB{https://www.lmfdb.org/ModularForm/GL2/Q/holomorphic/94/2/a/a/}{94.2.a.a}
& 
$d=10 = 2 \cdot 5 $
& 
-
&
$2 \leq n \leq 94$
&
\\ \hline
\multirow{2}{*}{96} 
& 
\LMRefDB{https://www.lmfdb.org/ModularForm/GL2/Q/holomorphic/96/2/a/a/}{96.2.a.a}
& 
$
\begin{matrix}
\eta_{96}[ 0, 0, 0, -2, 0, 4, 4, 0, -2, 0, 0, 0 ]\\ 
\eta_{96}[ 0, 0, 0, 4, 0, -2, -2, 0, 4, 0, 0, 0 ]
\end{matrix}
$
& 
$
\begin{matrix}
1 \\ -1
\end{matrix}
$
&
\checkmark
&
\\ 
\cline{2-6}
& 
\LMRefDB{https://www.lmfdb.org/ModularForm/GL2/Q/holomorphic/96/2/a/b/}{96.2.a.b}
& 
$
\begin{matrix}
\eta_{96}[ 0, 0, 0, -2, 0, 4, 4, 0, -2, 0, 0, 0 ]\\ 
\eta_{96}[ 0, 0, 0, 4, 0, -2, -2, 0, 4, 0, 0, 0 ]
\end{matrix}
$
& 
$
\begin{matrix}
1 \\ 1
\end{matrix}
$
&
\checkmark 
&
\\

\hline
98 
&
\LMRefDB{https://www.lmfdb.org/ModularForm/GL2/Q/holomorphic/98/2/a/a/}{98.2.a.a}
& 
$
\begin{matrix}
\eta_{98}[1,0,1,1,0,1] \\ 
\eta_{98}[0,1,1,1,1,0] \\ 
\eta_{98}[1,1,1,1,0,0] \\ 
\eta_{98}[0,0,1,1,1,1]
\end{matrix}
$
&
$
\begin{matrix}
4 \\ 4 \\ 1 \\ 7
\end{matrix}
$
& 
\checkmark
&
 \\ \hline 
\multirow{2}{*}{99} 
& 
\LMRefDB{https://www.lmfdb.org/ModularForm/GL2/Q/holomorphic/99/2/a/a/}{99.2.a.a}
& 
$d=2$
& 
-
&
$3\leq n \leq 4$
&
\\ 
\cline{2-6}
& 
\LMRefDB{https://www.lmfdb.org/ModularForm/GL2/Q/holomorphic/99/2/a/b/}{99.2.a.b}
& 
$d = 2$
& 
-
&
$3 \leq n \leq 18$
&
\\
\cline{2-6}
& 
\LMRefDB{https://www.lmfdb.org/ModularForm/GL2/Q/holomorphic/99/2/a/c/}{99.2.a.c}
& 
$
\begin{matrix}
\eta_{99}[-1, 4, -1, -1 ,4, -1] \\
\eta_{99}[1, 0, 1, 1, 0, 1] \\
\eta_{99}[0,2,0,0,2,0]
\end{matrix}$
& 
$\begin{matrix}
1 \\ 3 \\ -2
\end{matrix}$
&
\checkmark
&
\cite{AyginThesis}
\\
\cline{2-6}
& 
\LMRefDB{https://www.lmfdb.org/ModularForm/GL2/Q/holomorphic/99/2/a/d/}{99.2.a.d}
& 
$d = 2$
& 
-
&
$2 \leq n \leq 19$
&
\\

\hline
100 
&
\LMRefDB{https://www.lmfdb.org/ModularForm/GL2/Q/holomorphic/100/2/a/a/}{100.2.a.a}
& 
$
\begin{matrix}
\eta_{100}[ 0, 0, 0, 0, 2, 0, 0, 2, 0 ] \\ 
\eta_{100}[ 0, 2, 0, 0, 2, 0, 0, 0, 0 ] \\ 
\eta_{100}[ 0, 1, 0, 0, 2, 0, 0, 1, 0 ] 
\end{matrix}
$
&
$
\begin{matrix}
5 \\ 1 \\ 4
\end{matrix}
$
& 
\checkmark
&
 \\ \hline 
\end{tabular}
\end{table}

\end{document}